\documentstyle[12pt,leqno]{article}

\input boxedeps.tex
\SetepsfEPSFSpecial
\HideDisplacementBoxes
\def\reals{\hbox{\sl I\kern-.18em R \kern-.3em}}
\def\ints{\hbox{\sl Z\kern-.4em Z \kern-.3em}}
\def\nats{\hbox{\sl I\kern-.16em N \kern.05em}}
\def\rats{\hbox{\sl Q \kern-.83em\vrule height.59em
depth0em
\kern.87em}}
\def\complexes{\hbox{\sl\kern.50em I\kern-.50em C
\kern.05em}}

\newtheorem{theorem}{Theorem} 

\newtheorem{lemma}{Lemma}
\newtheorem{corollary}{Corollary}
\newtheorem{remark}{Remark}

\def\bs{\\ \bigskip}
\def\pf {{\bf Proof:} \ }
\def\endpf{$\|$ \bigskip}
\def\cK{{\cal K}}
\def\cF{{\cal F}}
\def\cT{{\cal T}}
\def\cTK{{\cal TK}}

\def\cS{{\cal S}}

\def\cM{{\cal M}}

\begin{document}
\title{On transversally simple knots}
\author {Joan S. Birman \thanks{The first author
acknowledges partial support from the U.S.National
Science Foundation under  Grants DMS-9705019 and
DMS-9973232.  The second author is a graduate student
in the Mathematics Department of Columbia University.
She was partially supported under the same grants.}  
 \and Nancy C. Wrinkle }
\date{March 1, 2001}

\maketitle
\centerline{Journal of Differential Geometry, to appear}

\flushleft

\begin{abstract}  This paper studies knots that are
transversal to the standard contact structure in
$\reals^3$, bringing techniques from topological knot
theory to bear on their transversal classification. We
say that a transversal knot type
$\cTK$ is {\it transversally simple} if it is 
determined by its topological knot type
$\cK$ and its Bennequin number. The main theorem
asserts that any
$\cTK$ whose associated $\cK$ satisfies a condition
that we call {\em exchange reducibility} is
transversally simple.

As a first application, we prove that the unlink is
transversally simple, extending the main theorem in
\cite{El}.  As a second application we use a new
theorem of Menasco \cite{Me} to extend a result of
Etnyre \cite{Et} to prove that all iterated torus knots
are transversally simple.  We also give a formula for
their maximum Bennequin number.  We show that the
concept of exchange reducibility is the simplest of the
constraints that one can place on $\cK$ in order to
prove that any associated $\cTK$ is transversally
simple.  We also give examples of pairs of transversal
knots that we conjecture are {\em not} transversally
simple.
\end{abstract}

\section{Introduction.}
\label{section:introduction} Let $\xi$ be the {\it
standard} contact structure in oriented 3-space
$\reals^3 = (\rho, \theta, z)$, that is the kernel of
$\alpha = \rho^2 d\theta + dz$.  An oriented knot
$K$ in contact $\reals^3$ is said to be a {\it
transversal} knot  if it is transversal to the planes
of this contact structure.  In this  paper, the term
`transversal' refers to this contact structure only. If
the knot
$K$ is parametrized by $(\rho(t),\theta(t), z(t))$,
then $K$ is transversal if and only if
$\frac {z'(t)}{\theta'(t)} \neq -(\rho(t))^2$ for every
t.  We will assume  throughout that
$\alpha > 0$ for all $t$,  pointing out later how our
arguments adapt to the case
$\alpha < 0$. \bs

\begin{figure}[htpb]
\centerline{\BoxedEPSF{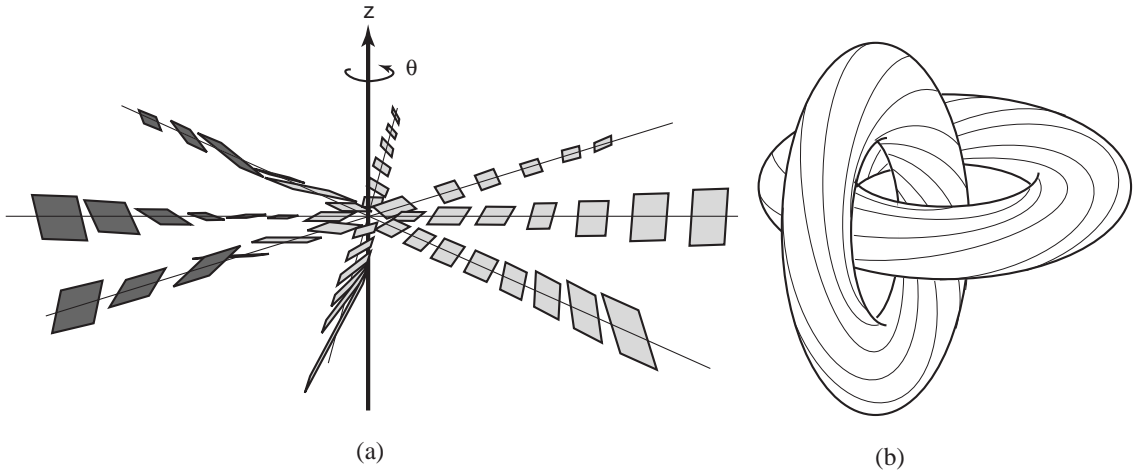 scaled 800}}
\caption {The standard contact structure on $\reals^3$
and the Hopf fibration on $S^3$.}
\label{figure:standard contact structure}
\end{figure}

For the benefit of the reader who may be unfamiliar
with the standard contact structure, Figure 
\ref{figure:standard contact structure}(a) illustrates
typical 2-planes in this structure in $\reals^3$, when
$z$ is fixed, and
$\rho$ and
$\theta$ vary.  The structure is radially symmetric. It
is also invariant under translation of $\reals^3$
parallel to the $z$ axis.  Typical 2-planes are
horizontal at points on the $z$ axis and twist
clockwise (if the point of view is out towards
increasing
$\rho$ from the $z$ axis,) as $\rho\to\infty$. \bs

There has been some discussion about whether the planes
tend to vertical as
$\rho\to\infty$ or to horizontal.  If one looks at the
limit of
$\alpha$, it appears that the limit is a rotation of
$\pi/2$.  However, if one derives the standard contact
structure on $S^3$ from the Hopf fibration, as
described below, and wants to have this structure be
consistent with the one defined on $\reals^3$, it is
necessary to take the limit to be a rotation up to (but
not through, as a rotation of more than $\pi$ results
in an overtwisted structure)
$\pi$.  The resulting contact structures on $\reals^3$
are equivalent, through a contactomorphism that
untwists the planes from $\pi$ to
$\pi/2$.  Thus we can work in the standard contact
structure on
$S^3$, which has horizontal planes in the limit, while
using the contact form $\alpha$, which induces vertical
planes in the limit.  The details are below. 

The standard contact structure extends to $S^3$ and has
an interesting interpretation in terms of the geometry
of $S^3$. Let 
$$S^3 = \{(z_1,z_2) = (\rho_1e^{i\theta_1},
\rho_2e^{i\theta_2})\in
\complexes^2 \ / \ \rho_1^2 + \rho_2^2 = constant.\}$$ 
Then $\xi$ is the kernel of $\rho_1^2d\theta_1 +
\rho_2^2d\theta_2$. The field of 2-planes may be
thought of as the field of hyperplanes which are
orthogonal to the fibers of the Hopf fibration
$\pi:S^3\to S^2$.  See Figure \ref{figure:standard
contact structure}(b) for a picture of typical fibers. 
Identify the $z$ axis in $\reals^3$ with the core of
one of the solid tori. There is a fiber through each
point in $S^3$, and the 2-plane at a point is
orthogonal to the fiber through the point. 
\bs

The (topological) {\it type} $\cK$ of a knot
$K\subset\reals^3$ is its equivalence class under
isotopy of the pair
$(K,\reals^3)$.  A sharper notion  of equivalence is
its {\it transversal knot type}
$\cTK$, which requires that $\frac {z'(t)}{\theta'(t)}
+ (\rho(t))^2$ be positive at every point of the
deformed knot during every stage of the isotopy.  The
difference between these two concepts is the central
problem studied in this paper. \bs

\begin{figure}[htpb]
\centerline{\BoxedEPSF{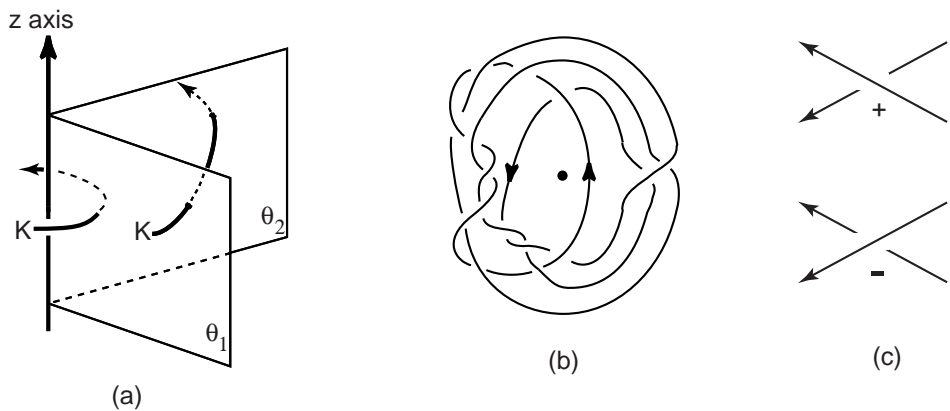 scaled 800}}
\caption {(a) Closed braid;  (b) Example of a closed
braid projection;   (c) Positive and negative crossings}
\label{figure:closed braid}
\end{figure}

A parametrized knot $K\subset\reals^3$ is said to be
represented as a {\em closed braid} if 
$\rho(t)>0$ and
$\theta^\prime(t)>0$ for all $t$.  See Figure
\ref{figure:closed braid}(a). It was proved by
Bennequin in $\S$23 of 
\cite{Be} that every transversal knot is transversally
isotopic to a transversal closed braid.   This result
allows us to apply results obtained in the study of
closed braid representatives of topological knots to
the problem of understanding transversal isotopy. We
carry Bennequin's approach one step further, initiating
a comparative study of the two equivalence relations:
topological equivalence of two closed braid
representatives of the same transversal knot type, via
closed braids, and transversal equivalence of the same
two closed braids. Transversal equivalence is of course
more restrictive than topological equivalence. \bs

Topological equivalence of closed braid representatives
of the same knot has been the subject of extensive
investigations by the first author and W. Menasco, who
wrote a series of six papers with the common title {\em
Studying links via closed braids}. See, for example,
\cite{BM3} and
\cite{BM5}. See also the related papers
\cite{BH} and
\cite{BF}. In this paper we will begin to apply what
was learned in the topological setting to the
transversal problem. See also Vassiliev's paper
\cite{V}, where we first learned that closed braid
representations of knots were very natural in analysis;
also our own contributions in
\cite{BW}, where we began to understand that there were
deep connections between the analytic and the
topological-algebraic approaches to knot theory. \bs

A well-known invariant of a transversal knot type
$\cT\cK$ is its {\it Bennequin} number $\beta(\cT\cK)$.
It is not an invariant of $\cK$.  We now define it in a
way that will allow us to compute it from a closed
braid representative $K$ of $\cT\cK$.  The {\em braid
index $n = n(K)$} of a closed braid $K$ is the linking
number of $K$ with the oriented $z$ axis. A generic
projection of $K$ onto the
$(\rho,\theta)$ plane will be referred to as a {\it
closed braid projection}. An example is given in Figure
\ref{figure:closed braid}(b). The origin in the
$(\rho,\theta)$ plane is indicated as a black dot; our
closed braid rotates about  the
$z$ axis in the direction of increasing
$\theta$. The {\it algebraic crossing number}
$e=e(K)$ of the closed braid is the sum of the signed
crossings in a closed braid projection, using the sign
conventions given in Figure
\ref{figure:closed braid}(c).   If the transversal knot
type
$\cTK$ is represented by a closed braid $K$, then its
Bennequin number
$\beta(\cTK)$ is:
$$ \beta(\cTK) = e(K) - n(K).$$ Since $e(K)-n(K)$ can
take on infinitely many different values as $K$ ranges
over the representatives of $\cK$, it follows that
there exist infinitely many transversal knot types for
each topological knot type. It was proved by Bennequin
in
\cite{Be} that $e(K)-n(K)$ is bounded above by
$-\chi(\cF)$, where
$\cF$ is a spanning surface of minimal genus for
$\cK$.  Fuchs and Tabachnikov gave a different upper
bound in
\cite{FT}.   However, sharp upper bounds are elusive
and are only known in a few very special cases. \bs

We now explain the geometric meaning of $\beta(\cTK)$.
Choose a point 
$(z_1,z_2) = (x_1+ix_2,x_3+ix_4)\in TK\subset S^3$.
Thinking of
$(z_1,z_2)$ as a point in $\reals^4$, let
$\vec{p} = (x_1,x_2,x_3,x_4)$. Let $\vec{q} =
(-x_2,x_1,-x_4,x_3)$ and let
$\vec{r} = (-x_3,x_4,x_1,-x_2)$. Then
$\vec{r}\cdot\vec{p} =
\vec{r}\cdot\vec{q}  = \vec{p}\cdot\vec{q}=0$.  Then
$\vec{q}$ may be interpreted as the outward-drawn
normal to the contact plane at
$\vec{p}$, so that $\vec{r}$ lies in the unique contact
plane at the point $\vec{p}\in S^3$.  Noting that a
transversal knot is nowhere tangent to the contact
plane, it follows that for each point
$\vec{p}$ on a transversal knot $TK\subset S^3$ the
vector $\vec{r}$ gives a well-defined direction for
pushing $TK$ off itself to a related simple closed
curve $TK'$. The Bennequin number
$\beta(\cTK)$ is the linking number
${\cal L}k(TK,TK')$. See $\S$16 of \cite{Be} for a
proof that
$\beta(TK)$ is invariant under transverse isotopy and
that $\beta(TK) = e(K)-n(K)$. \bs 

We say that a transversal knot is {\em transversally
simple} if it is characterized up to transversal
isotopy by its topological knot type and its Bennequin
number.  In
\cite{El} Eliashberg proved that a transversal unknot 
is transversally simple. More recently Etnyre
\cite{Et} used Eliashberg's techniques to prove that
transversal positive torus knots are transversally
simple.
\bs 

Our first main result, Theorem
\ref{theorem:realizing the maximal Bennequin number},
asserts that if a knot type $\cK$ is exchange
reducible, (a condition we define in Section 2,) then
its maximum Bennequin number is realized by any closed
braid representative of minimum braid index. As an
application, we are able to compute the maximum
Bennequin number for all iterated torus knots.  See
Corollary
\ref{corollary:maximum Bennequin numbers for iterated
torus knots}.  Our second main result, Theorem
\ref{theorem:exchange reducible implies t-simple}, 
asserts that if
$\cTK$ is a transversal knot type with  associated
topological knot type $\cK$, and if
$\cK$ is exchange reducible, then
$\cTK$ is transversally simple. As an application, we
prove in Corollary
\ref{corollary:iterated torus knots are transversally
simple} that transversal iterated torus knots are
transversally simple.  The two corollaries use new
results of Menasco
\cite{Me}, who proved (after an early version of this
paper was circulated) that iterated torus knots are
exchange reducible. In Theorem
\ref{theorem:3 braid examples} we establish the
existence of knot types that are {\em not} exchange
reducible.\bs

Here is an outline of the paper. Section
\ref{section:exchange reducibility and transversal
simplicity} contains our main results. In it we will
define the concept of an exchange reducible knot and
prove Theorems
\ref{theorem:realizing the maximal Bennequin number} and
\ref{theorem:exchange reducible implies t-simple}. In
Section
\ref{section:examples} we discuss examples,
applications and possible generalizations.
\bs

{\bf Acknowlegements} We thank Oliver Dasbach, William
Menasco and John Etnyre for conversations and helpful
suggestions relating to the work in this paper. We are
especially grateful to Menasco for the manuscript
\cite{Me}. In an early version of this paper, we
conjectured that iterated torus knots might be exchange
reducible. We explained our conjecture to him, and a
few days later he had a proof!  We also thank Wlodek
Kuperberg for sharing his beautiful sketch of the Hopf
fibration (Figure 1(b)) with us. 
Finally, we thank William Gibson, who noticed
our formula for the Bennequin number of iterated torus knots in an earlier
version of this paper and pointed out to us in a private conversation that it
could be related to the upper bound which was  given by Bennequin in \cite{Be},
by the formula  in Corollary 3, part (2).

\subsection{Remarks on techniques.} This subsection
contains a discussion of the techniques used in the
manuscripts \cite{BM4}, \cite{BM5}, \cite{BM94} and
\cite{Me}, tools which form the foundation on which the
results of this paper are based.  We compare these
techniques to those used in the manuscripts \cite{Be},
\cite{El} and \cite{Et}, although Bennequin's paper
rightfully belongs in both sets.  This description and
comparison is of great interest, but is not essential
for the reading of this paper or the digestion of its
arguments.  \bs  

We concern ourselves with two known foliations of an
orientable surface
$\cF$ associated to $K$: the {\it characteristic}
foliation $\xi_F$ from contact geometry and the {\it
topological} foliation from braid theory.  The
characteristic foliation of $\cF$ is the line field
$\xi \cap T\cF$, given by the intersection of the
planes of the contact structure with the planes of the
tangent space of the surface, which is then integrated
to a singular foliation of $\cF$.   The {\it
topological  foliation} is the foliation of $\cF$ which
is induced by intersecting the foliation of $\reals^3$
minus the $z$ axis (see Figure \ref{figure:the braid
structure},) with the surface
$\cF$.  The foliation of $\reals^3$ minus the $z$-axis
by half-planes is called the {\em standard braid
structure} on $\reals^3$ $-$ the
$z$-axis.  This structure is given by half-planes with
boundaries on the $z$ axis. The surface used in
\cite{El} and
\cite{BM5} was a spanning surface for $K$; in
\cite{BM4} it was a 2-sphere which intersects $K$
twice; in
\cite{BM94} it was a torus in the complement of $K$; in
\cite{Et} and
\cite{Me} it is a torus $T\subset S^3$ on which $K$ is
embedded.  Menasco also considers the foliation of a
meridian disc in the solid torus which $T$ bounds.  The
characteristic foliation of a surface (associated to a
transversal knot or to another transversal or
Legendrian curve,) is a tool of study in contact
geometry.  It was the main tool in the manuscripts
\cite{Be}, \cite{El} and
\cite{Et}.   \bs

In topological knot theory, one studies the topological
foliation of
$\cF$ defined above.  The review article \cite{BF} may
be useful to the reader who is unfamiliar with this
area.  The study of the topological  foliations has
produced many results, for example the classification
of knots that are closed 3-braids \cite{BM3} and a
recognition algorithm for the unknot \cite{BH}. Braid
theory was also a major tool in the work in \cite{Be},
but it appears that Bennequin's detailed study of the
foliation is based entirely on the characteristic
foliation, as it occurs for knots in $\reals^3$ and
$S^3$.  To the best of our knowledge this paper
contains the first application of braid foliations to
the study of  transversal knots.
\bs

\begin{figure}[htpb]
\centerline{\BoxedEPSF{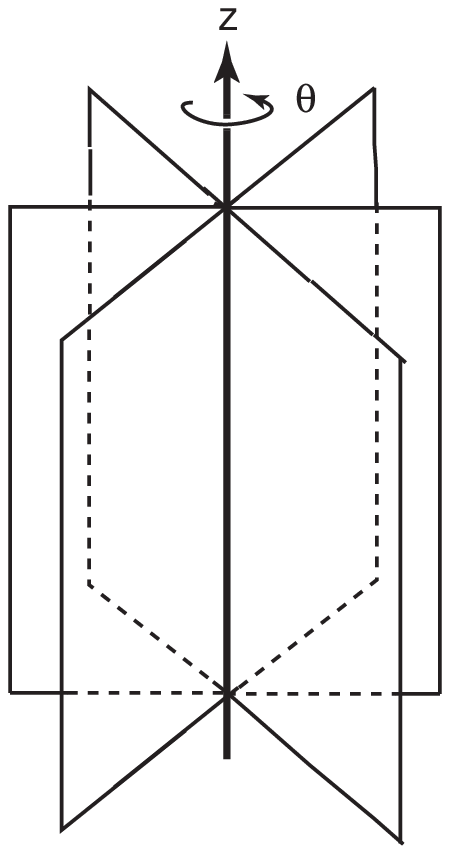 scaled 700}}
\caption {Half-planes in the braid structure on
$\reals^3$.}
\label{figure:the braid structure}
\end{figure}

We note some similarities between the two foliations:
The characteristic foliation is oriented and the braid
foliation is orientable.  (The orientation is ignored,
but a dual orientation, determined by an associated
flow,  plays an equivalent role.) The  foliations can
be made to agree in the limiting case, as
$\rho\to\infty$ (see the comments above Figure
\ref{figure:standard contact structure}). 
\bs

\begin{figure}[htpb]
\centerline{\BoxedEPSF{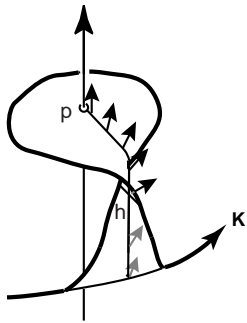 scaled 800}}
\caption {A positive elliptic singularity (p) and a
positive hyperbolic singularity (h).}
\label{figure:singularity sign}
\end{figure}

After an appropriate isotopy of
$\cF$ both foliations have no leaves that are simple
closed curves. Also, their  singularities are finite in
number, each being either an elliptic point or a
hyperbolic point (the hyperbolic point corresponding to
a saddle-point tangency of
$\cF$ with the 2-planes of the structure).   The signs
of the singularities of each foliation are determined
by identical considerations: the surface is naturally
oriented by the assigned orientation on the knot. If at
a singularity the orientation of the surface agrees
(resp. disagrees) with the orientation of the
foliation, then the singularity is positive (resp.
negative). See Figure
\ref{figure:singularity sign}.  \bs

In both foliations the hyperbolic singularities are
4-pronged singularities. If $s$ is a hyperbolic
singular point, then the four branches of the singular
leaf through $s$ end at either elliptic points or at a
point   on $K$.  (The condition that no singular leaves
of the characteristic foliation connect hyperbolic
points is a genericness assumption appearing in the
literature on Legendrian and transversal knots).   The
three possible cases are illustrated in Figure
\ref{figure:the three types of  hyperbolic
singularities}. In that figure the elliptic points are
depicted as circles  surrounding
$\pm$ signs (the sign of the elliptic singularity) and
the hyperbolic singularities  are depicted as black
dots.  Two of the four branches of the singular leaf
end at positive  elliptic points. The other two end at
either two negative elliptic points, or one negative
elliptic point and one point on
$K$, or two points on $K$. \bs

\begin{figure}[htpb]
\centerline{\BoxedEPSF{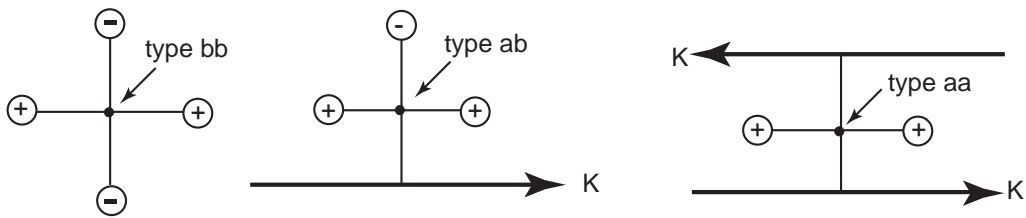 scaled 800}}
\caption {The three types of hyperbolic singularities.}
\label{figure:the three types of hyperbolic
singularities}
\end{figure}

There are also differences between the two foliations. 
In the braid foliation, elliptic points always
correspond to punctures of the surface by the
$z$ axis.  In the characteristic foliation, elliptic
points on the surface may or may not correspond to
punctures by the $z$-axis.  That is, there may be elliptic
points not corresponding to punctures, and there might be
punctures not corresponding to elliptic points.  Here is an
example.  In the braid foliation, if there is a piece of the
surface along the boundary, foliated by a single positive pair
of elliptic and hyperbolic singularities, then the only
possible embedding for that piece is shown in Figure
\ref{figure:singularity sign}.  On the other hand, in
the characteristic foliation, if there is a piece of
the surface, also along the boundary, also foliated by
a positive elliptic-hyperbolic pair, then the
corresponding embedding may or may not be the one shown
in Figure \ref{figure:singularity sign}.  The
embeddings will coincide if the tangent to the surface
at the $z$ axis is horizontal.  \bs In work on the
braid foliation one uses certain properties that appear
to have been ignored in work based upon the
characteristic foliation. For example, the work on
braid foliations makes much of the distinction between
the three types of hyperbolic singularities which we
just illustrated in Figure
\ref{figure:the three types of hyperbolic
singularities}, calling them {\it types bb, ab} and
{\it aa}.  The resulting combinatorics play a major
role in the study of braid foliations. It seems to us
that the distinction between $bb, ab$ and $aa$
singularities can also be made in the situation of the
characteristic foliation, but that this has not been
done. 
\bs

In the braid foliation the elliptic points have a
natural cyclic order on the $z$ axis, if we are considering the
ambient space as $S^3$ and the braid axis as one of the
core circles of the Hopf fibration, and the hyperbolic points
have a natural cyclic order in
$0\leq\theta\leq 2\pi$.   These orderings do not seem
useful in the contact setting.    On the other hand,
the characteristic foliation is invariant under
rotation by $\theta$ and translation by $z$, so the
interesting parameter seems to be the coordinate $\rho$.
\bs

An essential tool in manipulating and simplifying the
characteristic foliation is the Giroux Elimination
Lemma (\cite{Gi}, \cite{El}), which allows one to
`cancel' pairs of same sign singularities.  In
topological knot theory different modifications have
been introduced that are the braid foliation analogue
of isotopies of the Giroux Elimination Lemma, see
\cite{BM5} and also \cite{BF}. They are called  {\em
$ab$ exchange moves} and {\em
$bb$ exchange moves}, and they use pairs of Giroux-like
cancellations, but on a much larger scale.
\bs

\section{Exchange reducibility and transversal
simplicity.}
\label{section:exchange reducibility and transversal
simplicity}

Our initial goal is to motivate and define the concept
of exchange reducibility.  Let
$\cK$ be a topological knot type and let  $K$ be a
closed
$n$-braid representative of $\cK$. We consider the
following three modifications of $K$

\begin{itemize} 
\item  Our first modification is {\em braid isotopy},
that is, an isotopy in the complement of the braid
axis.  In \cite{Mo} it is proved that isotopy classes 
of closed
$n$-braids are in one-to-one correspondence with
conjugacy classes in the braid group
$B_n$.  Since the conjugacy problem in the braid group 
is a solved problem, each conjugacy class can then be
replaced by a unique representative that can be assumed
to be transversal.  Braid isotopy preserves the
Bennequin number since it preserves both braid index
and algebraic crossing number.

\item Our second move is {\em destabilization}. See
Figure 
\ref{figure:destabilization and exchange moves}(a). The
box labeled $P$ contains an arbitrary
$(n-1)$-braid, and the label $n-2$ on the braid strand
denotes $n-2$ parallel braid strands. The
destabilization move reduces braid index from
$n$ to
$n-1$ by removing a `trivial loop'.  If the trivial
loop contains a positive crossing, the move  is called
a {\em positive} or $+$ destabilization.  Positive
destabilization reduces algebraic crossing number and
preserves the Bennequin number.  Negative ($-$)
destabilization increases the Bennequin number by 2.

\item  Our third move is the {\em exchange move}. See
Figure
\ref{figure:destabilization and exchange moves}(b).  
\begin{figure}[htpb]
\centerline{\BoxedEPSF{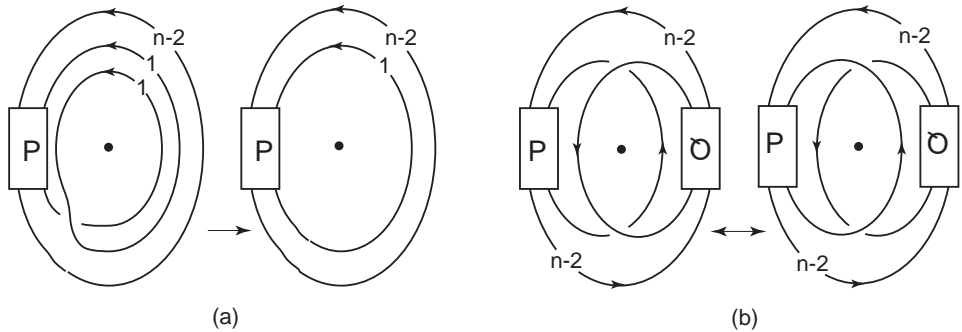 scaled 700}}
\caption {(a) positive destabilization and (b) The
exchange move}
\label{figure:destabilization and exchange moves}
\end{figure} In general the exchange move changes
conjugacy class and so cannot be replaced by braid
isotopy. The exchange move preserves both braid index
and algebraic crossing number, hence preserves the
Bennequin number. \bs
\end{itemize}

To motivate our definition of exchange reducibility, we
recall the following theorem, proved by the first
author and W. Menasco:
\bs

{\bf Theorem A} (\cite{BM5}, with a simplified proof in
\cite{BF}):  {\em Let $K$ be a closed
$n$-braid representative of the $m$-component unlink. 
Then
$K$ may be simplified to the trivial $m$-braid
representative, i.e. a union of $m$ disjoint round
planar circles, by a finite sequence of the following
three changes: braid isotopies, positive and negative
destabilizations, and exchange moves.} \bs

Motivated by Theorem A, we introduce the following
definition:\bs
\underline{ Definition:} A knot type $\cK$ is said to
be {\em exchange reducible} if an arbitrary closed
braid representative $K$ of arbitrary braid index
$n$ can be changed to an arbitrary closed braid
representative of minimum braid index
$n_{min}(\cK)$ by a finite sequence of  braid
isotopies, exchange moves and
$\pm$-destabilizations.  Note that this implies that
any two minimal braid index representatives are either
identical or are {\em exchange-equivalent}, i.e., are
related by a finite  sequence of braid isotopies and
exchange moves. \bs

Our first result is: 
\begin{theorem} 
\label{theorem:realizing the maximal Bennequin number}
If $\cK$ is an exchange reducible knot type, then the
maximum Bennequin number of $\cK$ is realized by any
closed braid representative of braid index
$n_{min}(\cK)$.
\end{theorem}

The proof of Theorem \ref{theorem:realizing the maximal
Bennequin number}
 begins with a lemma.  In what follows, we  understand
``transversal isotopy" to mean a topological isotopy
that preserves the condition 
$\alpha = \rho^2 d\theta + dz > 0$  at every point of
the knot and at every stage of the isotopy.  

\begin{lemma}
\label{lemma:transversal isotopies} If a transversal
closed braid is modified by one of the following
isotopies,  then the isotopy can be replaced by a
transversal isotopy: 

{\rm (1)} Braid isotopy.

{\rm (2)} Positive stabilization or positive
destabilization.

{\rm (3)} An exchange move.
\end{lemma} 

{\bf Proof of Lemma \ref{lemma:transversal isotopies}:}
\bs

Proof of (1): Since the braid strands involved in the
isotopy will be
$\gg
\epsilon$ away from the z-axis at each stage (so
avoiding $-\rho^2 = 0$), any isotopy will be
transversal if we keep the strands involved "relatively
flat" ($dz/d\theta \sim 0$) at each stage.  Since
everything is happening locally there is space to
flatten the strands involved without changing the
braid.\bs

\begin{figure}[htpb]
\centerline{\BoxedEPSF{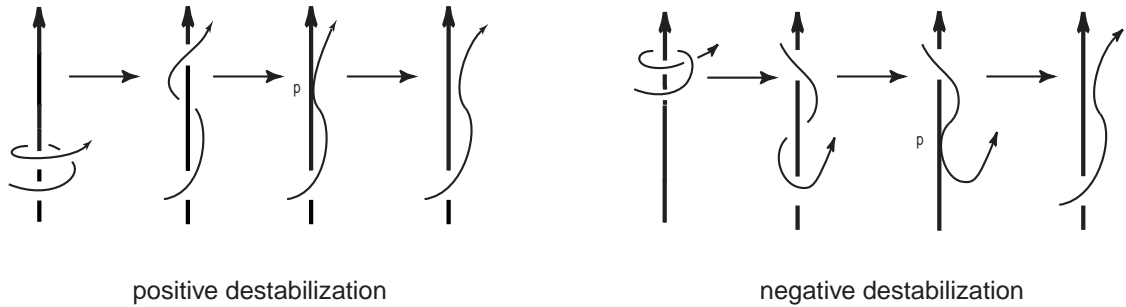 scaled 800}}
\caption {Destabilization, with a singularity at s,
where $d\theta = 0$ and $\rho = 0$.}
\label{figure:destabilization}
\end{figure}

Proof of (2):   See Figure
\ref{figure:destabilization}(a).  Consider a single
trivial loop around the $z$ axis, with a positive
crossing.  
 We have
$d\theta > 0$ along the entire length of the loop since
we are working with a closed braid.  For a positive
crossing we have $dz
\geq 0$ throughout the loop as well.  Therefore the
inequality
$dz/d\theta > - \rho ^2$ is true for all non-zero real
values of
$\rho$.  Crossing the $z$ axis to destabilize the braid
results in at least one singular point, where
$d\theta = 0$, but if we continue to keep $dz
\geq 0$ then in the limit, as $-\rho ^2
\rightarrow 0$ from the negative real numbers,
$dz/d\theta$ goes to 
$\infty$ through the positives.  Therefore
$dz/d\theta \neq -\rho ^2$ at any stage in the
isotopy.    \bs

\begin{figure}[htpb]
\centerline{\BoxedEPSF{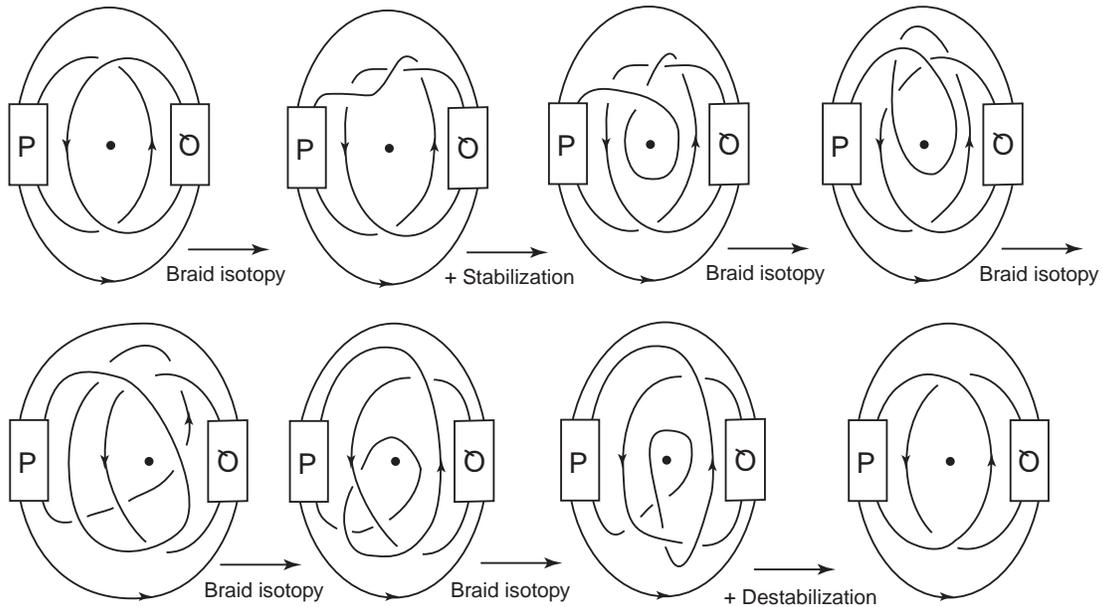 scaled 700}}
\caption {An exchange move corresponds to a sequence
consisting of braid isotopies,  a single positive
stabilization and a single positive destabilization.}
\label{figure:exchange moves}
\end{figure}

Proof of (3): The sequence of pictures in Figure
\ref{figure:exchange moves} shows that an exchange move
can be replaced by a sequence of the following moves:
isotopy in the complement of the
$z$ axis, positive stabilization, isotopy again, and
finally positive destabilization. Claim 3 then follows
from Claims 1 and 2.  \endpf  

\begin{remark}: {\rm Observe that the argument given to
prove (2) simply doesn't work for negative
destabilization.  See  figure
\ref{figure:destabilization}(b).   The singularity in
this destabilization is a point at which
$dz/d\theta = -\rho ^2$ in the limit. Indeed, a negative
destabilization can't be modified to one which is
transversal, because the Bennequin number (an invariant
of transversal knot type)  changes under negative
destabilization.} 
\end{remark}

\begin{remark} {\rm If we had chosen $\alpha < 0$ along
the knot, we would consider {\em negative}
stabilizations and destabilizations as transversal
isotopies and would use those instead of positive
stabilizations and destabilizations in the  exchange
sequence.  All the other moves translate to the
negative setting without change.}
\end{remark}

{\bf Proof of Theorem \ref{theorem:realizing the
maximal Bennequin number}:} Let $K$ be an arbitrary
closed braid representative of the  exchange reducible
knot type $\cK$.  Let
$K_0$ be a minimum braid index representative of
$\cK$, obtained from $K$  by the sequence described in
the definition of exchange reducibility.  We must prove
that the transversal knot
$TK_0$ associated to $K_0$ has  maximum Bennequin
number for the knot type
$\cK$. Note that in general $K_0$ is not unique,
however it will not matter, for if $K_0'$ is a
different closed braid representative of minimal braid
index, then $K_0$ and $K_0'$ are related by a sequence
of braid isotopies and exchange moves, both of which
preserve both braid index and algebraic crossing
number, so $\beta(K_0') =
\beta(K_0)$. 
\bs

We obtain $K_0$ from
$K$ by a sequence of  braid isotopies, exchange moves,
and
$\pm$-destabilizations.  Braid isotopy, exchange moves
and positive destabilization  preserve
$\beta(TK)$, but negative destabilization increases the
Bennequin number by 2,  so the sequence taking $K$ to
$K_0$ changes the Bennequin number  from
$\beta(TK) = c$ to
$\beta(TK_0) = c + 2p$, where p is the number of
negative destabilizations in the  sequence.  The
question then is whether $c + 2p$ is maximal for the
knot type $\cK$.  If
$c + 2p$ is less than maximal, then  there exists some
other closed braid representative
$K^\prime$ of the knot type $\cK$ with maximum 
$\beta(TK^\prime)  > \beta(TK_0)$.  Since $K_0$ has
minimum braid index for the knot type
$\cK$, it must be that $n(K^\prime) \geq n(K_0)$.   If
$n(K^\prime) = n(K_0)$, then the two braids are
equivalent by a sequence of Bennequin number preserving
exchange equivalences, so suppose instead that
$n(K^\prime) > n(K_0)$.  Then, since $K^\prime$ is a
closed braid representative of the exchange reducible
knot type
$\cK$, there must exist a sequence of braid isotopies,
exchange moves,  and
$\pm$-destabilizations taking $K^\prime$ to a minimum
braid index braid representative
$K_0 ^\prime$.  (As above, we also take any $K_0 ^\prime
\in
\cM_0$ ).  Since
$\beta(TK^\prime)$ is assumed to be maximum, and none
of the moves in  the sequence taking
$K^\prime$ to $K_0 ^\prime$ reduce Bennequin number, it
must be that
$\beta(TK_0 ^\prime) =
\beta(TK^\prime)$.  But since $K_0 ^\prime$ and
$K_0$ are both minimum braid  index representatives of
$\cK$, they must be equivalent by a sequence of
Bennequin number preserving exchange moves and 
isotopies.  Thus
$\beta(TK_0)=
\beta(TK_0^\prime)$. \endpf

Our next result is:

\begin{theorem}
\label{theorem:exchange reducible implies t-simple}  If
$\cTK$ is a transversal knot type with associated
topological knot type
$\cK$, where $\cK$ is exchange reducible, then
$\cTK$ is transversally simple.
\end{theorem}

The proof of Theorem \ref{theorem:exchange reducible
implies t-simple} begins  with two lemmas.   Our first
lemma had been noticed long ago by the first author and
Menasco, who  have had a long collaboration on the
study of closed braid representatives of knots and
links.  However, it had never been used in any of their
papers. It is therefore new to this paper,  although we
are indebted to  Menasco for his part in its
formulation.  A contact-theory analogue of Lemma
\ref{lemma:Oliver's lemma}, below, appears as Lemma 3.8
of \cite{Et}.

\begin{lemma}
\label{lemma:exchange moves are the obstruction}  Using
exchange moves and  isotopy in the complement of the
braid axis, one may slide a trivial loop on a closed
braid from one location to another on the braid.
\end{lemma}

\begin{figure}[htpb]
\centerline{\BoxedEPSF{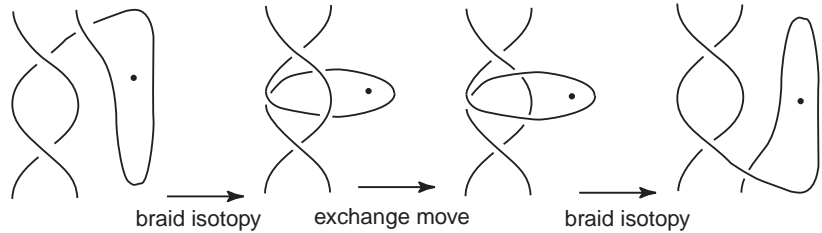 scaled 700}}
\caption {An exchange move that allows a negative
trivial loop to slide along a braid.}
\label{figure:aaexch}
\end{figure}

{\bf Proof of Lemma \ref{lemma:exchange moves are the
obstruction}:}   See Figure
\ref{figure:aaexch}.  It shows that, using braid
isotopy and exchange moves, we can slide a trivial
negative loop past any crossing to any  place we wish
on the braid. The argument for sliding a positive
trivial loop around  the braid is identical. 
\endpf
 
\begin{lemma}
\label{lemma:Oliver's lemma} Let $K_1$ and $K_2$ be
closed
$n$-braids that are exchange-equivalent.  Let
$L_1$ and $L_2$ be $(n+1)$-braids that are obtained from
$K_1$ and $K_2$ by either negative stabilization on
both or positive stabilization on both.  Then
$L_1$ and $L_2$  are exchange-equivalent.
\end{lemma} 

{\bf Proof of Lemma \ref{lemma:Oliver's lemma}:} We
already know there is a way to deform $K_1$ to
$K_2$, using  exchange equivalence.  Each braid isotopy
may be broken up into a  sequence of isotopies, each of
which only involves local changes on some well-defined
part of  the braid.  (For example, the defining
relations in the braid group are appropriate local 
moves on cyclic braids).  Similarly, exchange moves
have local support.  It may happen that the trivial
loop which we added interferes with the support of  one
of the isotopies or exchange moves. If so, then by
Lemma 
\ref{lemma:exchange moves are the obstruction} we may
use exchange equivalence to slide it out of the way. It
follows that we may deform
$L_1$ to $L_2$ by exchange equivalence. \endpf 

{\bf Proof of Theorem \ref{theorem:exchange reducible
implies t-simple}:}  We are given an arbitrary
representative of the transversal knot type $\cTK$. Let
$\cK$ be the  associated topological knot type.  By the
transversal Alexander's theorem
\cite{Be} we may  modify our representative
transversally to a transversal closed
$n$-braid $K = TK$ that represents the transversal knot
type $\cTK$ and the  topological knot type
$\cK$.  By the definition of exchange reducibility, we
may then find a  finite sequence of closed braids 
$$K = K_1\to K_2\to\dots \to K_{m-1}\to K_m,$$ all
representing
$\cK$, such that each $K_{i+1}$ is obtained from
$K_i$ by braid isotopy,  a positive or negative
destabilization or an exchange move, and such that
$K_m$ is a  representative of minimum braid index
$n_{min} = n_{min}(\cK)$ for the knot type
$\cK$.  The knots 
$K_1,\dots,K_m$ in the sequence will all have the
topological knot type
$\cK$. \bs

In general $\cK$ will have more than one closed braid
representative of  minimum braid index. Let
$\cM_0(\cK)=\{M_{0,1},M_{0,2}\dots\}$ be the set of all
minimum braid  index representatives of
$\cK$, up to braid isotopy. Clearly
$K_m\in\cM_0$.  {By \cite{Be}, each
$M_{0,i}\in\cM_0$ may be assumed to be a transversal
closed braid.
\bs

By Theorem \ref{theorem:realizing the maximal Bennequin
number} each
$M_{0,i}$ has maximal Bennequin number for all knots
that represent
$\cK$. In general this Bennequin number will not be the
same as the Bennequin  number of the original
transversal knot type $\cTK$.   By Lemma 
\ref{lemma:transversal isotopies} the moves that relate
any two 
$M_{0,i}, M_{0,j}\in \cM_0$ may be assumed to be
transversal. After all these modifications the closed
braids  in the set $\cM_0$ will be characterized, up to
braid isotopy, by their topological knot type
$\cK$, their braid index $n_{min}(\cK)$ and their
Bennequin number
$\beta_{max}(\cK)$. \bs  

If the transversal knot type $\cTK$ had Bennequin number
$\beta_{max}(\cK)$, it would necessarily follow that
$\cTK$ is characterized up to transversal isotopy by
its ordinary knot type and its Bennequin number. Thus
we have proved the theorem in the special case of
transversal knots that have maximum Bennequin number.
\bs 

We next define new sets $\cM_1,\cM_2,\dots$ of
transversal knots, inductively.  Each
$\cM_s$ is a collection of conjugacy classes of closed
$(n_{min}(\cK)+s)$-braids. We assume,  inductively,
that the braids in
$\cM_s$ all have topological knot type
$\cK$, braid index
$n_{min}(\cK)+s$ and Bennequin number 
$\beta_{max}(\cK)-2s$.  Also, their  conjugacy classes
differ at most by exchange moves. Also, the collection
of conjugacy classes of
$(n_{min}(\cK)+s)$-braids in the set
$\cM_s$ is completely determined by the collection of
conjugacy classes of  braids in the set
$\cM_0$.  We now define the set
$\cM_{s+1}$ by choosing an arbitrary closed braid
$M_{i,s}$ in
$\cM_s$ and  adding a trivial negative loop.  Of
course, there is no unique way to do this,  but by Lemma
\ref{lemma:exchange moves are the obstruction} we can
choose one such trivial  loop and use exchange moves to
slide it completely around the closed braid
$M_{i,s}$. Each  time we use the exchange move of Lemma
\ref{lemma:exchange moves are the obstruction}, we will
obtain a new conjugacy  class, which we then add to the
collection $\cM_{s+1}$. The set
$\cM_{s+1}$ is defined to be the collection of all
conjugacy classes of closed braids  obtained by adding
trivial loops in every possible way to each
$M_{i,s}\in\cM_s$. The closed braids  in
$\cM_{s+1}$ are equivalent under braid isotopy and
exchange moves. They all have topological knot type
$\cK$, braid index $n_{min}(\cK)+s+1$, and Bennequin
number 
$\beta_{max}(\cK) - 2(s + 1)$. The collection of closed
braids in the set
$\cM_{s+1}$ is completely determined by the collection
of closed braids in
$\cM_s$, and so by the closed braids in
$\cM_0$.   
\bs

In general negative destabilizations will occur in the
chain $K_1\to K_m$. Our plan is to change the order of
the moves in the sequence
$K_1\to K_m$, pushing all the negative destabilizations
to the right until we obtain a new sequence,  made up
of  two subsequences: 
$$K = K_1^\star\to K_2^\star\to\dots\to K_r^\star=
K_0^\prime\to\dots
\to K_s^\prime =  K_p,$$ where $K_p$ has minimum braid
index
$n_{min}(\cK)$. The first  subsequence 
$\cS_1$, will be 
$K = K_1^\star\to K_2^\star\to\dots\to K_r^\star$,
where every
$K_i^\star$ is a transversal representative of
$\cTK$ and the connecting moves are  braid isotopy,
positive destabilizations and exchange moves. The
second  subsequence, $\cS_2$, is
$K_r^\star = K_0^\prime\to\dots \to K_q^\prime$, where
every
$K_{i+1}^\prime$ is obtained from
$K_i^\prime$ by braid isotopy and a single negative
destabilization.  Also,
$K_q^\prime$  has minimum braid index
$n_{min}(\cK)$. 
\bs

To achieve the modified sequence, assume that 
$K_i\to K_{i+1}$ is the first negative destabilization.
If the negative trivial loop does not interfere with
the moves leading from
$K_{i+1}$ to $K_m$,  just renumber terms so that the
negative destabilization becomes
$K_m$ and every other
$K_j, j>i$ becomes $K_{j-1}$. But if it does interfere,
we need to slide it out of the way to remove the
obstruction. We use Lemma \ref {lemma:exchange moves
are the obstruction} to do that, adding exchange moves
as required. \bs  

So we may assume that we have our two subsequences
$\cS_1$ and
$\cS_2$. The braids in $\cS_1$ are all transversally
isotopic and so they all have the same Bennequin number
and they all represent
$\cTK$.   The braids $K^\prime_i
\in \cS_2$ all have the same knot type, but 
$\beta(K^\prime_{i+1}) = \beta(K^\prime_i) +2$ for each
$i = 1,
\dots, s-1$.  With each  negative destabilization and 
braid isotopy, the Bennequin number increases by 2 and
the braid index decreases by 1.   Each braid represents
the same knot type
$\cK$ but a different transversal knot type
$\cT\cK$. \bs 

Our concern now is with $\cS_2$, i.e. 
$K_r^\star = K_0^\prime\to\dots \to K_s^\prime$,  where
every
$K_{i+1}^\prime$ is obtained from $K_i^\prime$ by a
single negative  destabilization and braid isotopy. 
The number of negative  destabilizations in subsequence 
$\cS_2$ is exactly one-half the difference between the
Bennequin number
$\beta(\cTK)$ of the original transversal knot
$\cTK$ and the Bennequin number
$\beta_{max}(\cK)$. \bs

Let us now fix on any particular minimum braid index
representative of
$\cK$ as a minimum braid index closed braid
representative of  the transversal knot type that
realizes
$\beta_{max}(\cK)$. It will not matter which we choose,
because all belong  to the set
$\cM_0$ and so are exchange-equivalent. We may then
take the final braid
$K_r^\star$  in $\cS_1$, which is the same as the
initial braid
$K_0^\prime$  in $\cS_2$, as our representative of
$\cTK$, because it realizes the minimal  braid index
for $\cTK$  and by our construction, any other such
representative is related to the one we have chosen by
transversal isotopy. We  may also proceed back up the
sequence
$\cS_2$ from 
$K_s^\prime$ to a new representative that is obtained
from
$K_0^\prime$ by   adding $s$  negative trivial loops,
one at a time.  By repeated application of Lemma
\ref{lemma:Oliver's lemma} we know that choosing any
other element of
$\cM_0$ will take us to an exchange-equivalent element
of $\cM_s$.  In this way we arrive in the set
$\cM_s$,  which also contains $K_r^\star$, and which is
characterized by $\cK$ and
$\beta$. The proof of Theorem
\ref{theorem:exchange reducible implies t-simple} is
complete. 
\endpf

\section{Examples, applications and possible
generalizations.}
\label{section:examples} In this section we discuss
examples which illustrate Theorems 
\ref{theorem:realizing the maximal Bennequin number} 
and
\ref{theorem:exchange reducible implies t-simple}.

\subsection{The unlink and the unknot.} Theorem A,
quoted earlier in this manuscript, asserts that the
$m$-component unlink, for $m \geq 1$, is exchange
reducible. In considering  a link transversally, it
should be mentioned that we are assuming each of the
components of the link satisfy the same inequality
$\alpha > 0$.  We also need to define the Bennequin
number properly for this transversal link.  The natural
way to do so, suggested by Oliver Dasbach, is by the
following method.  For a crossing involving two
different components of the link, assign
$\pm 1/2$ to each component depending on the sign of
the crossing.  Assign $\pm 1$ to each crossing
consisting of strands from the same component, as in
the case of a knot.  Then the Bennequin number of each
component is the difference between the algebraic
crossing number
$e$ (a sum of $\pm 1$'s and
$\pm 1/2$'s) and $n$, the braid index of that
component.  Define the Bennequin number of the link to
be the collection of the Bennequin numbers of the
components of the link.  The following corollary is an
immediate consequence of Theorems
\ref{theorem:realizing the maximal Bennequin number}
and A:

\begin{corollary}
\label{corollary:the unlink is transversally simple}
The  $m$-component unlink, $m\geq 1$, is transversally
simple. In particular, the unknot is transversally
simple.
\end{corollary} Note that Corollary
\ref{corollary:the unlink is transversally simple}
gives a new proof of a theorem of Eliashberg \cite{El}.

\subsection{Torus knots and iterated torus knots}  In the manuscript
\cite{Et} J. Etnyre proved that positive torus knots
are transversally simple.  His proof failed
for negative torus knots, but he conjectured that the
assertion was true for all torus knots and possibly
also for all iterated torus knots.  In an early draft
of this manuscript we conjectured that torus knots and
iterated torus knots ought to be exchange reducible,
and sketched our reasons.  Happily, the conjecture is
now a fact, established by W. Menasco in \cite{Me}.  
Two corollaries follow.  To state and prove our first
corollary, we need to fix our conventions for the description of torus
knots and iterated torus knots.
\bs

\underline{Definition:} Let $U$ be the unit circle in the plane $z=0$,
and let 
$N(U)$ be a solid torus of revolution with
$U$ as its core circle. Let
$\lambda_0$ be a longitude for $U$, i.e.
$\lambda_0$ is a circle in the
plane $z=0$ which lies on $\partial N(U)$, so
that $U$ and
$\lambda_0$ are concentric circles in the plane $z=0$.   See Figure
\ref{figure:the standard solid torus}. A {\em torus knot} of  {\em type}
$e(p,q)$, where $e=\pm$, on $\partial N(U)$, denoted
$K_{e(p,q)}$, is the closed $p$-braid
$(\sigma_1\sigma_2\cdots\sigma_{p-1})^{eq}$ on
$\partial N(U)$, where
$\sigma_1,\dots,\sigma_{p-1}$ are elementary braid
generators of the braid group $B_{p}$.  Note that
$K_{e(p,q)}$ intersects the curve $\lambda_0$ in
$q$ points, and note that the algebraic crossing
number of its natural closed braid projection on the
plane $z=0$ is $e(p-1)q$.  The knot
$K_{e(p,q)}$ also has a second natural closed braid
representation, with the unknotted circle $U$ as braid axis and the closed
$q$-braid $(\sigma_1\sigma_2\cdots\sigma_{q-1})^{ep}$ as closed
braid representative. Since $p$ and $q$ are coprime integers, one of these
closed braids will have smaller braid index than the other, and without loss
of generality we will assume in the pages which follow that we have
chosen $p$ to be smaller than $q$, so that $K_{e(p,q)}$ has braid index $p$.  

\begin{figure}[htpb]
\centerline{\BoxedEPSF{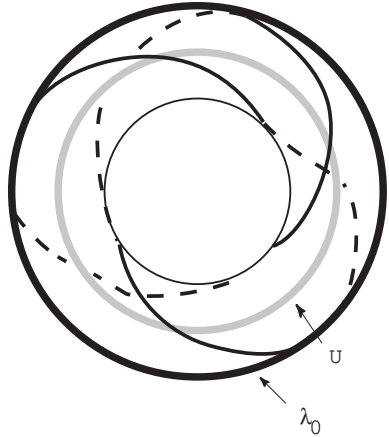 scaled 700}}
\caption {The standard solid torus $N(U)$, with
$K_{+(2,3)}\subset\partial N(K)$}
\label{figure:the standard solid torus}
\end{figure}

\underline{Definition}: We next define what we mean by
an {\em $e(s,t)$-cable} on a knot $X$ in 3-space. Let
$X$ be an arbitrary oriented knot in oriented
$S^3$, and let $N(X)$ be a solid torus neighborhood of
$X$ in 3-space. A {\em longitude} $\lambda$ for  $X$ is
a simple closed curve on
$\partial N(X)$ which is homologous to
$X$ in $N(X)$ and null-homologous in
$S^3\setminus X$. Let $f:N(U)\to N(X)$ be a
homeomorphism which maps
$\lambda_0$ to $\lambda$. Then
$f(K_{e(s,t)})$ is an $e(s,t)$-{\em cable} about
$X$.\bs

\underline{Definition}: Let $\{e_i(p_i,q_i), i =
1,\dots,r\}$ be a choice of signs
$e_i=\pm$ and coprime positive integers 
$(p_i,q_i)$, ordered so that for each
$i$ we have  $p_i,q_i>0$.   An {\it iterated torus knot}
$K(r)$ of type
$(e_1(p_1,q_1),\dots,e_r(p_r,q_r))$, is defined
inductively by:
\begin{itemize}
\item $K(1)$ a torus knot of type
$e_1(p_1,q_1)$, i.e. a type $e(p_1,q_1)$ cable on the unknot $U$.  Note that,
by our conventions, $p_1<q_1$.
\item $K(i)$ is an  $e_i(p_i,q_i)$ cable about
$K(i-1)$.  We place no restrictions on the relative magnitudes of
$p_i$ and $q_i$ when $i>1$.
\end{itemize}

Here is one of the simplest non-trivial examples of an
iterated torus knot.  Let $K(1)$ be the positive
trefoil, a torus knot of type $(2,3)$.  See
Figure
\ref{figure:trefoil as a torus knot}(a) and (b).  Note
that in the left sketch the core circle is our unit
circle $U$, while in the right sketch the core circle
is the knot $K(1)$.  The iterated torus knot
$K(2) = K_{(2,3),-(3,4))}$ is the $-(3,4)$ cable about
$K(1)$.  See Figure \ref{figure:(3,4) cabling}.  \bs
\begin{figure}[htpb]
\centerline{\BoxedEPSF{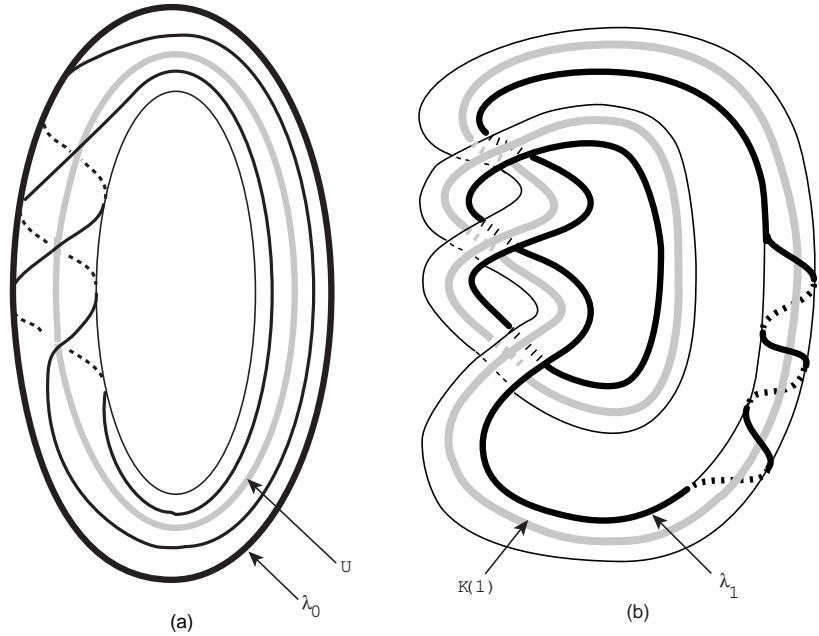 scaled 600}}
\caption {(a) the torus knot $K_{+(2,3)}$.  (b)the
solid torus neighborhood $N(K(1))$ of
$K(1)$, with core circle $K(1)$ and longitude
$\lambda_1$ marked.}
\label{figure:trefoil as a torus knot}
\end{figure} 
\begin{figure}[htpb!]
\centerline{\BoxedEPSF{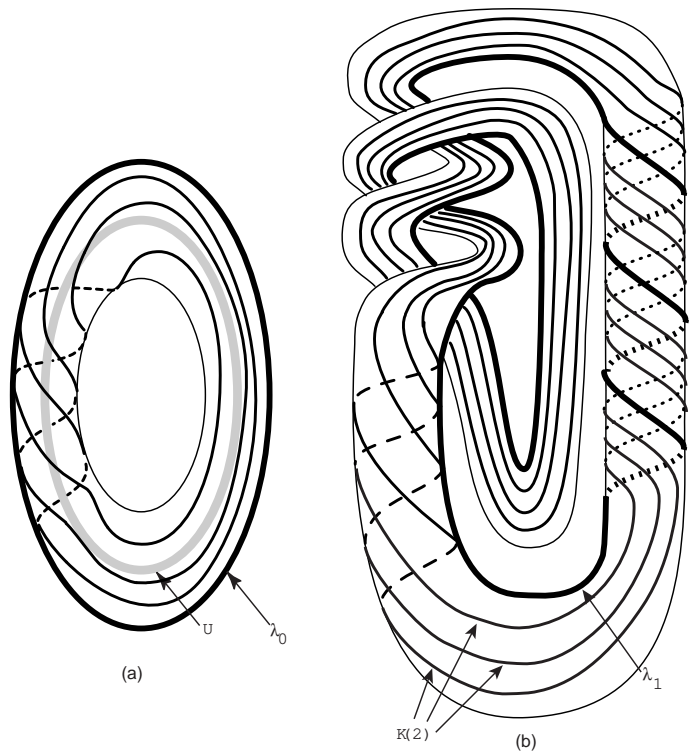 scaled 550}}
\caption {(a) the torus knot of type -(3,4) embedded in
$\partial N(U)$, (b) the iterated torus knot
$K(2) = K_{+(2,3), -(3,4)}$}
\label{figure:(3,4) cabling}
\end{figure}

In \cite{Me}, it was shown that the braid
foliation machinery used for the torus in \cite{BM94}
could be adapted to the situation in which the knot is
on the surface of the torus.  The main result of that
paper is the following theorem. \bs

{\bf Theorem:} (\cite{Me}, Theorem 1) {\it Oriented iterated torus
knots are exchange reducible.}\bs

Combining Menasco's Theorem with Theorem
\ref{theorem:exchange reducible implies t-simple}, we
have the following immediate corollary:

\begin{corollary}
\label{corollary:iterated torus knots are transversally
simple} Iterated torus knots are transversally simple.
\end{corollary}

Our next contribution to the theory of iterated torus knots requires that we
know the braid index of an iterated torus knot. The formula is implicit in the
work of Schubert \cite{Sch}, but does not appear explicitly there.

\begin{lemma}
\label{lemma:braid index of iterated torus knots}
Let $K(r) = K_{e_1(p_1,q_1),\dots,e_r(p_r,q_r)}$ be an r-times iterated torus
knot. Then the braid index of
$K(r)$ is $p_1p_2\cdots p_r$.
\end{lemma}

\pf We begin with the case $r=1$. By hypothesis $p_1<q_1$, also the torus
knot $K_{e_1(p_1,q_1)}$ is represented by a $p_1$-braid
$(\sigma_1\sigma_2\cdots\sigma_{p_1-1})^{e_1q_1}$.  By the formula given in
\cite{Jo} for the HOMFLY polynomial of torus knots, together with the
Morton-Franks-Williams braid index inequality (discussed in detail in
\cite{Jo}), it follows that this knot cannot be represented as a closed
$m$-braid for any $m<p_1$.  \bs

Passing to the general case, Theorem 21.5 of \cite{Sch} tells us that the
torus knot $K_{e_1(p_1,q_1)}$ and the array of integers
$e_1(p_2,q_2),\dots,e_r(p_r,q_r)$ form a complete system of invariants of the
iterated torus knot $K(r) = K_{e_1(p_1,q_1),\dots,e_r(p_r,q_r)}$.  Lemma
23.4 of \cite{Sch} tells us that, having chosen a $p_1$-braid representative
for $K_{e_1(p_1,q_1)}$, there is a natural $p_1p_2\cdots p_r$-braid
representative of $K(r)$. This representative is the only one on this
number of strings, up to isotopy in the complement of the braid axis. Theorem
23.1 of \cite{Sch} then asserts that $K(r)$ also cannot be represented as a
closed braid with fewer strands.  That is, its braid index is $p_1p_2\cdots
p_r$. \endpf

{\bf Remark:} The iterated torus knot $K_r$ has two natural closed braid
representatives. The first is a $p_1p_2\cdots p_r$- braid which has
the core circle $U'$ of the unknotted solid torus $S^3\setminus N(U)$ as braid
axis. The second is a $q_1p_2\cdots p_r$-braid which has the core circle
$U$ of the unknotted solid torus $N(U)$ as braid axis. In the case
$r=1$, the second choice gives a closed braid
which is reducible in braid index, i.e. it has $q_1-p_1$ trivial loops. From
this it follows that if $r>1$ it will have $(q_1-p_1)p_2\cdots p_r$ trivial
loops, thus the second closed braid representation is reducible to the first. 

\

We are now ready to state our second corollary about iterated torus knots. Let
$\chi$ be the Euler characteristic of an oriented surface of minimum genus
bounded by $K(r)$.

\begin{corollary}
\label{corollary:maximum Bennequin numbers for iterated
torus knots} 
Let $K(r) = K_{e_1(p_1,q_1),\dots,e_r(p_r,q_r)}$ be an iterated torus knot,
where $p_1<q_1$. Then the maximum Bennequin number of 
$K(r) = K_{e_1(p_1,q_1),\dots,e_r(p_r,q_r)}$ is given by the following two
equivalent formulas: 
\begin{enumerate}
\item [{\rm (1)}] $\beta_{max}(K(r)) = a_r-p_1p_2\cdots p_r,$  where
$a_r = \sum_{i=1}^r e_iq_i(p_i - 1)p_{i+1}p_{i+2}\ldots p_r.$
\item [{\rm (2)}] $\beta_{max}(K(r)) = -\chi - d$, 
where $d = \sum_{i=1}^r {(1-e_i)(p_i - 1)q_ip_{i+1}p_{i+2}
\ldots p_r}$.  
\item [{\rm (3)}] Moreover, the upper bound in the inequality $\beta_{max}(K(r))
\leq -\chi$ is achieved if and only if all of the $e_i's$ are positive.
\end{enumerate}
\end{corollary}

\

\pf We begin with the proof of (1). By Lemma \ref{lemma:braid index
of iterated torus knots}  the braid index of $K(r)$ is $p_1p_2\cdots p_r$. 
Therefore
$\beta_{max}(K(r)) = a_r - p_1p_2\cdots p_r$, where
$a_r$ is the algebraic crossing number of the unique $p_1p_2\cdots p_r$-braid
representative of $K(r)$.  To compute $a_r$ we proceed inductively.  If
$r=1$ then $K(1)$ is a type $e_1(p_1,q_1)$ torus knot, which is
represented by the closed  $p_1$-braid
$(\sigma_1\sigma_2\dots\sigma_{p-1})^{e_1q_1}$. Its
algebraic crossing number is $a_1 = e_1(p_1-1)q_1$. \bs

The knot $K(i)$ is an $e_i(p_i,q_i)$ cable on
$K(i-1)$.  Note that $K_{e_i(p_i,q_i)} \subset
\partial N(K_0)$, also 
$K_{e_i(p_i,q_i)}$ is a $p_i$-braid, also
$K_{e_i(p_i,q_i)}\cap\lambda_0$ consists of
$q_i$ points.  We shall think of the projection of the
braided solid torus
$N(K(i-1))$, which is a $p_1p_2\cdots p_{i-1}$-braid,
as being divided into three parts. The reader may find
it helpful to consult Figures
\ref{figure:iterated torus knots} (a), (b), (c) as we
examine the contributions to $a_i$ from each part.

\begin{figure}[htpb]
\centerline{\BoxedEPSF{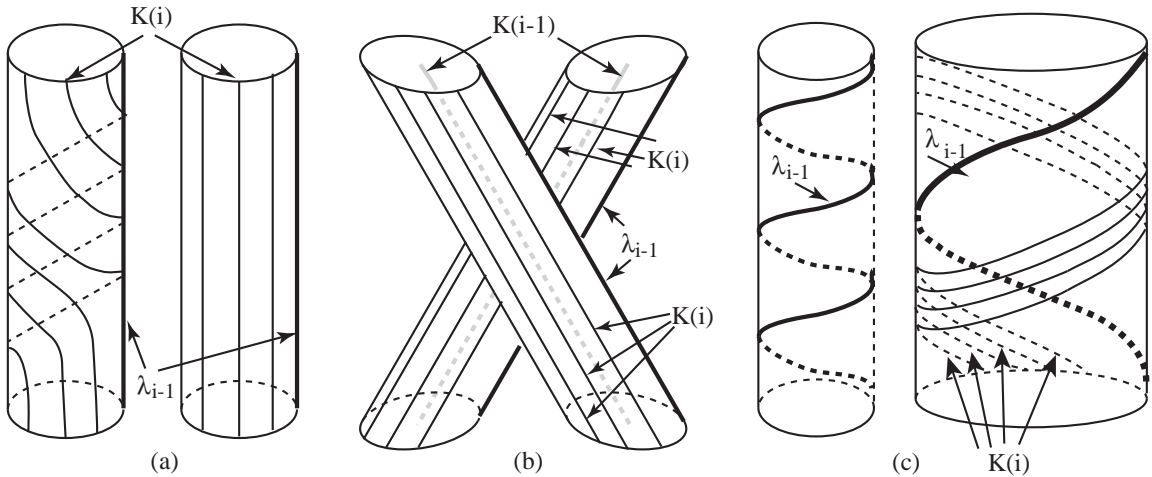 scaled 800}}
\caption {Iterated torus knots.}
\label{figure:iterated torus knots}
\end{figure}

\begin{itemize}
\item [(a)] The first part of $N(K(i-1))$ is the
trivial $p_1p_2\cdots p_{i-1}$-braid. The longitude
$\lambda_{i-1}$ is parallel to the core circle of
$N(K(i-1))$ in this part. The surface $\partial
N(K(i-1))$ contains on one of its $p_1p_2\cdots
p_{i-1}$ cylindrical branches the image under $f$ of
the braided part of
$K_{e_r(p_r,q_r)}$. See Figure 
\ref{figure:iterated torus knots}(a), which shows the
braid when
$e_i(p_i,q_i) = -(3,4)$. This part of
$K(i)$ contributes 
$e_i(p_i-1)q_i$ to $a_i$.   Note that there are
$q_i$ points where $f(K_{e_i(p_i,q_i)})$ intersects
$\lambda_{i-1}.$ 

\item [(b)] The second part of $N(K(i-1))$ contains all
of the braiding in $K(i-1)$, and so also in
$N(K(i-1))$. In Figure
\ref{figure:iterated torus knots}(b) we have
illustrated a single crossing in $K(i-1)$ and the
associated segments of $N(K(i-1))$. We show a single
crossing of
$\lambda_{i-1}$ (as a thick line) over
$K(i-1)$.  The single signed crossing  contributes
$p_i^2$ crossings to
$K(i)$, so the total contribution from all of the
crossings in
$K(i-1)$ will be $a_{i-1}p_i^2$.  The illustration
shows the case
$p_i=3$. 

\item [(c)] The third part of $N(K(i-1))$ is again the
trivial 
$p_1p_2\cdots p_{i-1}$-braid. It contains  corrections
to the linking number of
$\lambda_{i-1}$ with $K(i-1)$ which result from the
fact that a curve which is everywhere `parallel' to the
core circle will have linking number $a_{i-1}$, not 0,
with
$K(i-1)$. To correct for this, we must allow the
projected image of
$\lambda_{i-1}$ to loop around $\partial N(K(i-1))$
exactly
$-a_{i-1}$ times, so that its total linking number with
$K(i-1)$ is zero.  See the left sketch in Figure
\ref{figure:iterated torus knots}(c), which shows the 3
positive loops which occur if $a_{i-1}=-3$.
We have already introduced $q_i$ intersections between
$f(K_{e_i(p_i,q_i)})$ and $\lambda_{i-1}$, and
therefore we must avoid any additional intersections
which might arise from the
$-a_{i-1}$  loops. See the right sketch in Figure
\ref{figure:iterated torus knots}(c).  When
$\lambda_{i-1}$ wraps around $N(K(i-1))$ the additional
$-a_{i-1}$ times, the
$p_i$-braid $f(K_{e_i(p_i,q_i)}$ must follow. Each loop
in
$\lambda_{i-1}$ introduces $(p_i-1) + (p_i-2) + \cdots
+ 2 + 1 =
\frac{p_i(p_i-1)}{2}$ crossings per half-twist in
$f(K_{e_i(p_i,q_i)})$. Since there is a full twist to
go around the positive loop this number is doubled to
$p_i(p_i-1)$.  We have shown the 12 crossings in
$K(i)$ which come from a single loop in
$\lambda_{i-1}$ when $p_i=4$. The total contribution is
$-a_{i-1}p_i(p_i-1)$.
\end{itemize}

Adding up all these contributions we obtain
$$a_i = e_i(p_i-1)q_i +
a_{i-1}(p_i)^2-a_{i-1}(p_i^2 - p_i) = e_i(p_i-1)q_i +
a_{i-1}p_i.$$ Summing the various terms to compute $a_r$, we have proved part
(1) of the Corollary. 

\bigskip

The proof of (2) will follow from that of (1) if we can show that:
$$ \chi = p_1p_2\cdots p_r - d - a_r,$$
where $d = \sum_{i=1}^r {(1-e_i)(p_i - 1)q_ip_{i+1}p_{i+2}
\ldots p_r}$.  To see this, we must find a natural surface of minimum genus
bounded by $K(r)$ and compute its Euler characteristic.  By Theorem 12, Lemma
12.1 and Theorem 22 of \cite{Sch}, a surface of minimum genus bounded by $K(r)$
may be constructed by Seifert's algorithm, explained in Chapter 5 of
\cite{Ro}, from a representative of K(r) which has minimal braid index. We
constructed such a representative in our proof of Part (1). To compute its
Euler characteristic, use the fact that $\chi$ is the number of Seifert circles
minus the number of unsigned crossings (Exercises 2 and 10 on pages 119 and 121
of \cite{Ro}). By a theorem of Yamada
\cite{Ya} the number of Seifert circles is the same as the braid index, i.e.
$p_1p_2\cdots p_r$ in our situation.   The number of unsigned crossings is
$b_r$, where $b_i = (p_i-1)q_i + b_{i-1}p_i$ and
$b_1 = (p_1-1)q_1$ and $\chi = p_1p_2\cdots p_r - b_r$.   
Adding up the contributions from all the
$b_i's$ we get $b_r = \sum_{i=1}^r {(p_i - 1)q_ip_{i+1}p_{i+2}
\ldots p_r}$ which can be rewritten as $\sum_{i=1}^r
{(1-e_i+e_i)(p_i - 1)q_ip_{i+1}p_{i+2} \ldots p_r}$. 
Separating terms: 

$$b_r =  
\sum_{i=1}^r {(1-e_i)(p_i - 1)q_ip_{i+1}p_{i+2} \ldots p_r}+
\sum_{i=1}^r {(e_i)(p_i - 1)q_ip_{i+1}p_{i+2} \ldots p_r} = d
+ a_r.$$   

The claimed formula for $\chi$ follows.\bs

To prove (3), observe that the only case when $\beta_{max} = -\chi$ exactly
occurs when $d=0$, i.e. the sum
$\sum_{i=1}^r {(1-e_i)(p_i - 1)q_ip_{i+1}p_{i+2} \ldots p_r}
= 0$, That is, all the
$e_i$'s are $+1$.  \endpf

\subsection{Knots that are not exchange reducible}
\label{subsection:knots that are not exchange
reducible} A very naive conjecture would be that all
knots are exchange reducible, however that is far from
the truth. We begin with a simple example.  In the
manuscript \cite{BM3} Birman and Menasco studied knots
that are represented by closed 3-braids, up to braid
isotopy, and identified the proper subset of those
knots whose minimum braid index is 3 (i.e. not 2 or 1).
They prove that the knots that have minimum braid index
representatives of braid index 3 fall into two groups:
those that have a unique such representative (up to
braid isotopy) and infinitely many examples that have
exactly two distinct representatives, the two being
related by 3-braid {\em flypes}.  A flype is the knot 
type preserving isotopy shown  in  Figure
\ref{figure:the 3-braid flype}.  Notice that the flype
is classified as  {\em positive} or {\em negative}
depending on the sign of the isolated crossing.   After
staring at the figure, it should become clear to the
reader that the closed braids in a `flype pair' have 
the same topological knot type.  We say that a braid
representative {\em admits a flype} if  it is conjugate
to a braid that has the special form illustrated in
Figure
\ref{figure:the 3-braid flype}. 

\begin{figure}[htpb]
\centerline{\BoxedEPSF{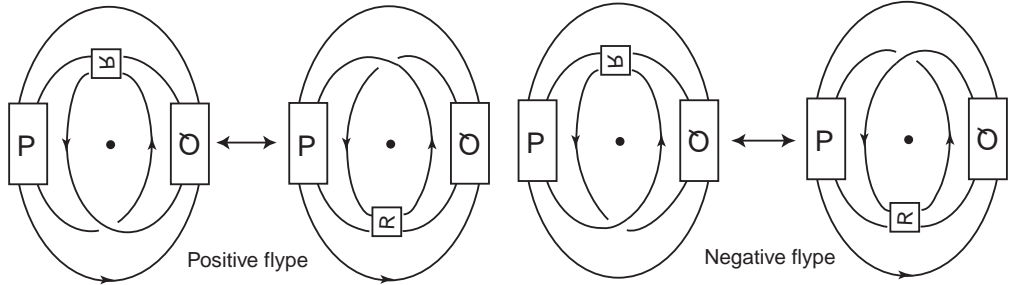 scaled 700}}
\caption {Closed 3-braids that are related by a flype.}
\label{figure:the 3-braid flype}
\end{figure}

\begin{theorem} {\rm \cite{BM3}}
\label{theorem:3 braid examples} The infinite sequence
of knots of braid index 3 in
{\rm \cite{BM3}}, each of which has two closed 3-braid
representatives, related by flypes, are examples of
knots that are {\em not} exchange reducible.
\end{theorem}

\pf It is proved in \cite{BM3} that for all but an
exceptional set of
$P,Q,R$ the  closed braids in a flype pair are in
distinct conjugacy classes. Assume from now on that a
`flype pair' means one of these non-exceptional pairs. 
Since conjugacy classes are in one-to-one
correspondence with braid isotopy equivalence classes,
it follows that the braids in a flype pair are not
related by braid isotopy.  On the other hand, it is
proved in
\cite{BM3} that when the braid index is
$\leq 3$ the exchange move can always be replaced by
braid isotopy, so the braids in a flype pair cannot be
exchange equivalent. \endpf

On closer inspection, it turns out that a positive
flype can be replaced by a sequence of braid isotopies
and positive stabilizations and destabilizations, which
shows that it is a transversal isotopy.  See Figure
\ref{figure:positive flypes are transversal}.
 
\begin{figure}[htpb]
\centerline{\BoxedEPSF{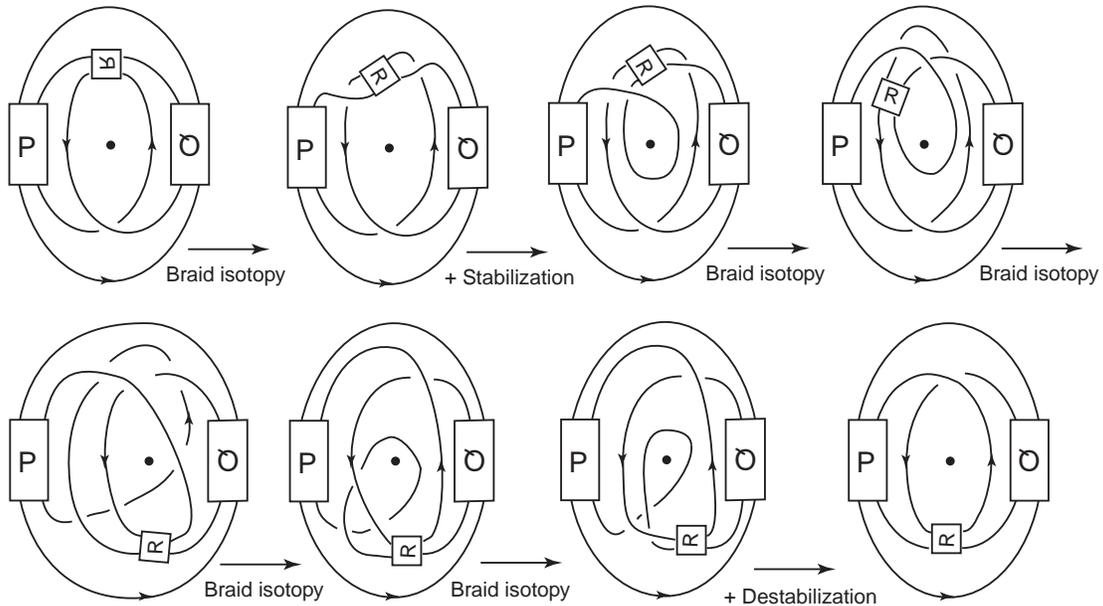 scaled 700}}
\caption {A positive flype can be replaced by a
sequence of transversal isotopies.}
\label{figure:positive flypes are transversal}
\end{figure}

The figure is a generalization of Figure
\ref{figure:exchange moves} (proving that exchange
moves are transversal), because flypes are
generalizations of exchange moves.  We replace one of
the
$\sigma^{\pm1}_n$  with the braid word we label
$R$.  A negative flype also has a replacement sequence
similar to the one pictured in Figure
\ref{figure:positive flypes are transversal},  but the
stabilizations and destabilizations required are
negative.    Therefore the  negative flype sequence
cannot be  replaced by a transversal isotopy using
these methods.  There may well be some other
transversal isotopy that can replace a negative flype,
but we did not find one.   Thus we are lead to the
following conjecture:
\bs

{\bf Conjecture:} {\em Any transversal knot type whose
associated topological knot type
$\cK$  has a minimum braid index representative that
admits a negative flype is {\it not} transversally
simple.} \bs

The simplest example which illustrates our conjecture
is shown  in Figure \ref{figure:flype example}. \bs
\begin{figure}[htpb]
\centerline{\BoxedEPSF{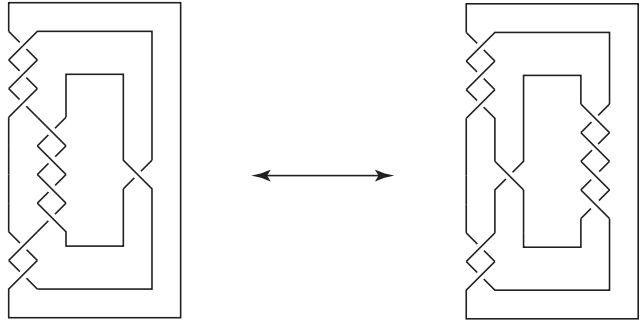 scaled 600}}
\caption {The simplest example which illustrates our
conjecture}
\label{figure:flype example}
\end{figure}

The essential difficulty we encountered in our attempts
to prove or disprove this conjecture is that the only
effective invariants of transversal knot type that are
known to us at this writing are the topological knot
type and the Bennequin number, but they do not
distinguish these examples.  \bs

\subsection{Knots with infinitely many transversally
equivalent closed braid representatives, all of minimal
braid index:}  At this writing the only known examples
of transversally simple knots are iterated torus knots.
By Theorem 24.4 of \cite{Sch} iterated torus knots have
unique closed braid representative of minimum braid
index, and it follows from this and Theorem
\ref{theorem:realizing the maximal Bennequin number}
that they have unique representatives of maximum
Bennequin number. It seems unlikely to us that all
transversally simple knots have unique closed braid
representatives of minimum braid index, and we now
explain our reasons.
\bs

The exchange move was defined in Figure
\ref{figure:destabilization and exchange moves}(b) of 
$\S$\ref{section:exchange reducibility and transversal
simplicity}. It seems quite harmless, being nothing
more than a special example of a Reidemeister move of
type II. It also seems  unlikely to produce infinitely
many examples of anything, however that is exactly what
happens when we combine it with braid isotopy. See
Figure
\ref{figure:infinitely many closed braids} with
$n = 4$. 
\begin{figure}[htpb]
\centerline{\BoxedEPSF{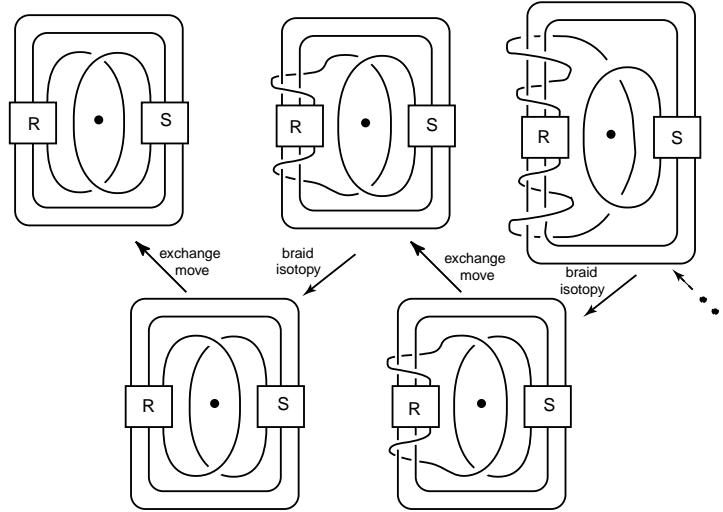 scaled 700}}
\caption {The exchange move and braid isotopy can lead
to infinitely many distinct closed
$n$-braid representatives of a single knot type. }
\label{figure:infinitely many closed braids}
\end{figure}
 Proceeding from the right to the left and following
the arrows, we see how braid isotopy and exchange moves
can be used to produce infinitely many examples of
closed braids which are transversally equivalent. It is
not difficult  to choose the braids $R$ and
$S$ in Figure \ref{figure:infinitely many closed
braids} so that the resulting closed braids are all
knots, and also so that they actually have braid index
4, and also so that they are in infinite many distinct
braid isotopy classes (using an invariant of Fiedler
\cite{Fi} to distinguish the braid isotopy classes). 
We omit details because, at this writing, we do not
know whether the knots in question are exchange
reducible, so we cannot say whether they all realize
the maximum Bennequin number for their associated knot
type.

\subsection{Generalizing the concept of `exchange
reducibility'}
\label{subsection:remarks on exchange reducibility}
Some remarks are in order on the concept of `exchange
reducibility'. Define two closed braids
$A\in B_n$, $A'\in B_m$  to be {\em Markov-equivalent}
if the knot types defined by the closed braids
coincide.  Markov's well-known theorem (see
\cite{Ma35}) asserts that
$Markov$-equivalence is generated by braid isotopy,
$\pm$-stabilization, and $\pm$-destabilization.
However, when studying this equivalence relation one
encounters the very difficult matter that
$\pm$-stabilization is sensitive to the exact spot on
the closed braid at which one attaches the trivial
loop.  On the other hand, Lemma
\ref{lemma:exchange moves are the obstruction} shows
that exchange moves are the obstruction to moving a
trivial loop from one spot on a knot to another.
Therefore if we allow exchange moves in addition to
braid isotopy, $\pm$-stabilization, and
$\pm$-destabilization, one might hope to avoid the need
for stabilization. That is the idea behind the
definition of exchange reducibility, and behind the
proof of Theorem A.  However, that hope is much too
naive, as was shown by the examples in
$\S$\ref{subsection:knots that are not exchange
reducible}.
\bs

A way to approach the problem of transversal knots is
to augment the definition of exchange reducible by
allowing additional `moves'. In their series of papers
{\em Studying knots via closed braids I-VI}, the first
author and Menasco have been working on generalizing
the main result in Theorem A to all knots and links. In
the forthcoming manuscript
\cite{BM01}, a general version of the `Markov theorem
without stabilization' is proved. The theorem states
that for each braid index
$n$ a finite set of new moves suffices to reduce any
closed braid representative of any knot or link to
minimum braid index `without stabilization'.  These
moves include not only exchange moves and positive and
negative flypes, but more generally {\em handle
moves} and {\em G-flypes}. Handle moves can always be realized
transversally. The simplest example of a G-flype is the 3-braid
flype that is pictured in Figure
\ref{figure:the 3-braid flype}, with weights assigned to
the strands. This will change it to an $m$-braid flype, for any
$m$.  But other examples exist, and they are much more
complicated.  We note, because it is relevant to the
discussion at hand, that any {\em positive} G-flype can be realized 
by a transversal isotopy. The sequence in Figure
\ref{figure:positive flypes are transversal} is a proof
of the simplest case. Awaiting the completion of 
\cite{BM01}, we leave these matters for future investigations.

\small {Joan S. Birman, Math. Dept, Barnard College of
Columbia University,}
\small {604 Mathematics, 2990 Broadway, New York, N.Y.
10027. e-mail: jb@math.columbia.edu}\bs

\small {Nancy C. Wrinkle, Math. Dept., Columbia
University,}
\small {408 Mathematics, 2990 Broadway, New York, N.Y.
10027. e-mail: wrinkle@math.columbia.edu}

\end{document}